\documentclass{amsart}

\usepackage{amsaddr}

\usepackage{amssymb}
\usepackage{amsmath}
\usepackage{color}
\usepackage{graphicx}
\usepackage{comment}
\usepackage{float}

\DeclareMathOperator*{\argmin}{\text{argmin}}

\def\E{{\mathbb E}}
\def\P{{\mathbb P}}
\def\R{{\mathbb R}}

\def\N{{\mathbb N}}

\newtheorem{theo}{Theorem}
\newtheorem{proposition}{\indent Proposition}

\title{Exponential inequality for chaos based on sampling without replacement.}

\vspace{1.5cm}

\author{P. Hodara and P. Reynaud-Bouret}
\address{INRA Jouy-en-Josas and Universit\'e C\^ote d'Azur, CNRS, LJAD, France}

\begin{document}

\maketitle

\begin{abstract}
We are interested in the behavior of particular functionals, in a framework where the only source of randomness is a sampling without replacement. 
More precisely the aim of this short note is to prove  an exponential concentration inequality for special U-statistics of order 2, that can be seen as chaos. 
\end{abstract}

\section{Introduction}
Since the introduction by \cite{Fisher} of permuted samples or by  Efron of bootstrapped procedures \cite{Efron, EfronTib},  many works have been devoted to the  study of the conditional  distribution of functionals of a given set of observations submitted to resampling procedures. Most of the time, asymptotic results are given (see for instance \cite{AG:89}). 

For  some special type of wild bootstrap,  one can go a bit further and provide easily non asymptotic concentration inequality because the resampling pattern is in some sense i.i.d. This leads to very powerful adaptive results, which need this non-asymptotic property to control several bootstrapped distribution at once. For instance, \cite{arlot} provides adaptive multiple tests and confidence regions in symmetric settings and \cite{FLRB} provides adaptive tests of equality for Poisson processes. However, in\cite{FLLR2012}, the authors have shown that in more straightforward set-ups, the correct bootstrap procedures that provide the right distribution are based on permutation of the data or sampling without replacement and this for U-statistics of order 2. The lack of concentration inequalities for such functionals prevented them to compute further properties such as separation rates.

If the asymptotic properties of empirical processes or U-statistics under random permutation or bootstrap procedures are quite well understood \cite{Hoeff51,AG:92}, advances to provide concentration inequalities are more rare. Permutation being intimately linked to exchangeability, there is a pioneer result in \cite{Aldous}, which expresses for empirical processes, the fact that sampling without replacement and sampling with replacement are linked and this link may be used to control one by another. Much more recently (see \cite{albert} and the references therein), it has been shown that conditional distribution of empirical processes under random permutation do have a Bernstein type concentration inequality.

Here we want to go one step beyond and provide exponential inequality for a special type of U-processes, of chaos type,  under sampling without replacement. This problem can also be seen as a particular functional of random permutation, where the only property that matters after permutation of $n$ observations is whether the new index is after $n/2$ or not (see Section 2.1 for more details). In this sense, it can also be seen as a first step for providing  concentration inequalities for U-statistics under random permutation. The main argument is based on a coupling between the position indicator and Rademacher variables, which has been used for the first time in \cite{CR:16} (see Section 3.1 for more details). This allows us to couple our statistics with a Rademacher chaos, for which concentration is known for many years, especially in the work of Latala (see for instance  \cite{latala} or the precursor results in \cite{dlPG}).

\section{Notation and main result}
Let $n$ be a positive even integer and $\left( a_{ij} \right)_{1 \leq i \neq j \leq n}$ a finite family of real numbers.

Pick at random without replacement exactly $p=n/2$ indices in $\{1,...,n\}$ and declare that 
\begin{itemize}
\item $\varepsilon_i=1$ if $i$ is picked, 
\item $\varepsilon_i=-1$ if it is not the case.
\end{itemize}

We are interested in the following chaos
\begin{equation}\label{chaosbizarre}
 Z= \sum_{i\neq j} \varepsilon_i\varepsilon_j a_{i,j}.
 \end{equation}

\subsection{Bootstrap point of view}

As an illustration from a bootstrap point of view, the $a_{i,j}$ might be think of as $g(X_i,X_j)$ for a certain bounded function $g$. In this set-up, we are given a full set of distinct observations $(X_1,...,X_n)$ that is divided in two samples $S_1=\{X_1,...,X_p\}$ and $S_2=\{X_{p+1},...,X_n\}$ for $n=2p$ and we are interested by a particular U-statistic of the form

$$U(X_1,...,X_n)=\sum_{\begin{array}{c} X_i \mbox{ and  } X_j \\\mbox{in the same sample} \end{array}} g(X_i,X_j) - \sum_{\begin{array}{c} X_i \mbox{ and  } X_j \\ \mbox{in different samples} \end{array}} g(X_i,X_j),$$ 

that is typical of two-sample problems \cite{FLLR2012}. If the two samples have the same distribution, then one can bootstrap the distribution of $U$ by drawing without replacement in $(X_1,...,X_n)$ a new sample $(X^*_1,...,X^*_n)$, which is in fact up to random permutation exactly the original sample. The unconditional distribution of $U$ is unaffected by the permutation and the main question is how the bootstrap, conditional distribution given the original set of observations $(X_1,...,X_n)$ behaves. Note that the conditional distribution of $U(X^*_1,...,X^*_n)$ given $(X_1,...,X_n)$ is exactly the one of $Z$ and that therefore exponential inequality for $Z$ gives upper bound on the conditional quantile of the bootstrap distribution.

\subsection{Rademacher chaos}
If we replace the $\varepsilon_i$'s in \eqref{chaosbizarre} by i.i.d. Rademacher variables $\varepsilon_i'$ that are equal to $1$ or $-1$ with probability $1/2$, we obtain a classical Rademacher chaos
$$Z'=\sum_{i\neq j} \varepsilon_i'\varepsilon_j' a_{i,j}.$$
 A typical exponential result for such a chaos is Corollary 3.2.6. in \cite{dlPG}, which states that there exists a constant $\kappa>0$ such that
$$ \E\left[ e^{\frac{|Z'|}{\kappa \sigma}}\right] \leq 2,$$
 with $\sigma^2= \sum_{i\neq j} a_{i,j}^2$.
 This leads to the following exponential inequality: for all $x>0$
\begin{equation}
\label{chaos-exp}
\P(|Z'| \geq \kappa \sigma  x) \leq 2 e^{-x}.
\end{equation}
Other refined results have been proved (see for instance \cite{latala}).

 \subsection{Main result}
Thanks to a coupling detailed  in the next section and due to \cite{CR:16}, we have proved the following result.
 
  \begin{theo}
  \label{main}
  There exists positive absolute constants $c$ and $C$ such that for all $x>0$
  $$\P\left( Z \geq c n  \max_{i\neq j} |a_{i,j}| (x+\log(n)) \right) \leq Ce^{-x}.$$
  \end{theo}
  
 Note that if for some positive $M$, the $a_{i,j}$'s are equal to $M$ or $-M$ at random, then up to constants, Theorem \ref{main} and \eqref{chaos-exp} coincide and in this sense Theorem \ref{main} gives the right order of magnitude for $x$ larger than $\log(n)$. However there is of course a significative loss  when $n^{-2}\sum_{i\neq j} a_{i,j}^2$ is much smaller than $\max_{i\neq j} |a_{i,j}|^2$.  A more refined inequality can be found 
 in Proposition \ref{gene}.

\section{Proof}

\subsection{Coupling}
The proof is based on a coupling argument, already used in \cite{CR:16}.
Let $(\varepsilon_i')_{i \in \N_*}$ be a sequence of i.i.d. Rademacher variables and let us define the stopping time $T$ by

\begin{equation}\label{eq:defineT}
T:=\argmin_t \left\{ \max \left( \sum_{i=1}^t 1_{\varepsilon_i'=1} \, , \, \sum_{i=1}^t 1_{\varepsilon_i'=-1} \right) = \frac{n}{2} \right\}
\end{equation}

We can interpret the Rademacher variables as labels giving whether the observation $X_i$ in a bootstrap paradigm, is put  in the first or second new sample. In this sense, a  i.i.d. Rademacher sampling means that we do not care about the fact that both new sample have the same size $p=n/2$. If we want to transform this sampling into another sampling of exactly the same size, we need to stop the i.i.d. Rademacher sampling at time $T$, give observation $X_T$ to the sample corresponding to $\varepsilon'_T$  and give all the remaining observations to the other sample corresponding to $-\varepsilon'_T$.

It means that we can define the $\varepsilon_i$ in the definition \eqref{chaosbizarre} by,
for each $1 \leq i \leq n,$ 
\begin{equation}\label{eq:defEpsPrime}
\varepsilon_i = \varepsilon_i' 1_{i \leq T} - \varepsilon_T' 1_{i>T}.
\end{equation}

Thanks to this coupling, one can prove the following Proposition.
\begin{proposition}\label{gene}
For every integer $0\leq \delta \leq n$, for all $x>0$,
\begin{multline*}
\P\left(|Z|\geq \kappa x \sqrt{\sum_{i\neq j \leq n-\delta} a_{i,j}^2}  + \sqrt{2 \left(x+\log(\delta) \right)} \left(\sum_{j > n-\delta} \sqrt{\sum_{i \leq n-\delta} a_{i,j}^2} + \sum_{i> n-\delta} \sqrt{\sum_{j \leq n-\delta} a_{i,j}^2} \right)+ \sum_{i\neq j > n-\delta}|a_{i,j}|\right) \\ \leq 2 e^{\frac{- \delta^2}{2(n-\delta)}}+6e^{-x},
\end{multline*}
where $\kappa$ is given by Corollary 3.2.6. of \cite{dlPG} (see also \eqref{chaos-exp}).
\end{proposition}

Note that as shown in the following proofs, it seems that the main term of the previous inequality is $$\sqrt{2 \left(x+\log(\delta) \right) } \left(\sum_{j> n-\delta} \sqrt{\sum_{i\leq n-\delta} a_{i,j}^2}+ \sum_{i>  n-\delta} \sqrt{\sum_{j\leq n-\delta} a_{i,j}^2} \right)$$ and not the one coming from the Rademacher chaos, i.e. $$\kappa x\sqrt{\sum_{i\neq j \leq n-\delta} a_{i,j}^2}  $$ at least for convenient choices of $\delta$, typically $\delta\simeq \sqrt{n x}$.  Even the last term $\sum_{i\neq j > n-\delta}|a_{i,j}|$ is not negligible with respect to the Rademacher chaos trend (see \eqref{chaos-exp}). We do not know if a more refined argument based on the same coupling may lead to a main term which is of the same order as \eqref{chaos-exp}, or if this loss is intrinsic to this coupling argument.

\subsection{Proof of Proposition \ref{gene}}

First of all,  let us control $\P(T<n-\delta)$.
$$
 \P \left( T \leq n-\delta \right) = \P \left( \max \left\{ \sum_{i=1}^{n-\delta} 1_{\varepsilon_i'=1} , \sum_{i=1}^{n-\delta} 1_{\varepsilon_i'=-1} \right\} \geq \frac{n}{2} \right) \leq 2 \P \left( \sum_{i=1}^{n-\delta} 1_{\varepsilon_i'=1} \geq \frac{n}{2} \right).
$$
But 
$$\P \left( \sum_{i=1}^{n-\delta} 1_{\varepsilon_i'=1} \geq \frac{n}{2} \right)=\P\left(\sum_{i=1}^{n-\delta} 1_{\varepsilon_i'=1} - \E(1_{\varepsilon_i'=1}) \geq \frac{\delta}{2}\right).$$
Hence using Hoeffding inequality, we obtain
\begin{equation}\label{eq:ctrlTsmall}
\P \left( T \leq n-\delta \right) \leq 2 e^{-\frac{\delta^2}{2(n-\delta)}}.
\end{equation}
Next on the event $\{T>n-\delta\}$, $\varepsilon_i=\varepsilon_i'$ for all $i\leq n-\delta$. This leads to, for any $u,v,w,z \in \R_+,$
\begin{multline}\label{eq:decomposition}
\P \left( \sum_{i\neq j \leq n} \varepsilon_i \varepsilon_j a_{ij} \geq u+v+w+z\right) \leq 
\P \left( T \leq n-\delta \right) + \P \left( \sum_{j\neq i \leq n-\delta} \varepsilon_i' \varepsilon_j' a_{ij} \geq u  \right) \nonumber \\ 
+ \P \left( \sum_{i \leq n- \delta < j} \varepsilon_i' \varepsilon_j a_{ij} \geq v \right) +  \P \left( \sum_{j \leq n- \delta < i} \varepsilon_i \varepsilon_j' a_{ij} \geq w \right) 
+  \P \left(\sum_{n- \delta < i\neq j} \varepsilon_i \varepsilon_j a_{ij} \geq z  \right). 
\end{multline}

The term $\sum_{j\neq i \leq n-\delta} \varepsilon_i' \varepsilon_j' a_{ij}$ is a classical Rademacher chaos and we can apply \eqref{chaos-exp} directly with $u=\kappa \sqrt{\sum_{i\neq j \leq n-\delta} a_{i,j}^2} x$ and
$\P \left( \sum_{j\neq i \leq n-\delta} \varepsilon_i' \varepsilon_j' a_{ij} \geq u  \right) \leq 2 e^{-x}.$

For the last term, since we are past $n-\delta$ and therefore roughly speaking past time $T$, it is likely that all the $\varepsilon_i$ are constant. Therefore, it is sufficient to take $z= \sum_{n- \delta < i\neq j} |a_{i,j}|$ and $\P \left(\sum_{n- \delta < i\neq j} \varepsilon_i \varepsilon_j a_{ij} \geq z  \right)=0$.

For $\P \left( \sum_{i \leq n- \delta < j} \varepsilon_i' \varepsilon_j a_{ij} \geq v \right)$, remark that
$$ \P \left( \sum_{i \leq n- \delta < j} \varepsilon_i' \varepsilon_j a_{ij} \geq v \right) \leq \P \left( \sum_{j>n-\delta} \left | \sum_{i \leq n- \delta } \varepsilon_i'  a_{ij} \right| \geq v \right).$$
By applying Hoeffding's inequality, we obtain that for every $j$,  and every $y>0$,
$$\P \left(  \left | \sum_{i \leq n- \delta } \varepsilon_i'  a_{ij} \right| \geq \sqrt{2y \sum_{i\leq n-\delta} a_{i,j}^2}\right)\leq 2 e^{-y}.$$
Therefore with
$$v= \sum_{j>n-\delta}  \sqrt{2y \sum_{i\leq n-\delta} a_{i,j}^2},$$
$$\P \left( \sum_{j>n-\delta} \left | \sum_{i \leq n- \delta } \varepsilon_i'  a_{ij} \right| \geq v \right)\leq 
\P\left(\exists j >n-\delta, \left | \sum_{i \leq n- \delta } \varepsilon_i'  a_{ij} \right| \geq \sqrt{2y \sum_{i\leq n-\delta} a_{i,j}^2}\right) \leq 2 \delta e^{-y}.$$
It remains to take $y=x+\log(\delta)$ and to proceed similarly for $\sum_{j \leq n- \delta < i} \varepsilon_i \varepsilon_j' a_{ij}$ to conclude the proof.

\section{Proof of Theorem \ref{main}}
We apply Proposition \ref{gene} and
first of all, we need to choose $\delta$. With $\delta=\sqrt{2nx}$, we have that $e^{-\frac{\delta^2}{2(n-\delta)}}\leq e^{-x}$.

Next, using $M=\max_{i\neq j} |a_{i,j}|$,  there exists an absolute  positive constant $c$ such that
\begin{multline*}
\kappa \sqrt{\sum_{i\neq j \leq n-\delta} a_{i,j}^2} x + \sqrt{2 \left( x+\log(\delta) \right) } \left(\sum_{j> n-\delta} \sqrt{\sum_{i\leq n-\delta} a_{i,j}^2} + \sum_{i> n-\delta} \sqrt{\sum_{j\leq n-\delta} a_{i,j}^2} \right)+ \sum_{i\neq j > n-\delta}|a_{i,j}| \\\leq 
c\left(n M x +  \delta M \sqrt{n}\sqrt{x+\log(nx)} +\delta^2 M\right)
\end{multline*}

It remains to replace $\delta$ by its $\sqrt{2 nx }$ to conclude.

\section*{Acknowledgements}
This research has been conducted as part of FAPESP project {\em Research, Innovation and
Dissemination Center for Neuromathematics} (grant 2013/07699-0 and grant 2016/17655-8).
This work was also supported by the French government, through the UCA$^{Jedi}$ "Investissements d'Avenir" managed by the National Research Agency (ANR-15-IDEX-01), by the structuring program $C@UCA$ of Universit\'e C\^ote d'Azur and by the interdisciplinary axis MTC-NSC of the University of Nice Sophia-Antipolis.

\bibliographystyle{plain}

\end{document}